\pdfoutput=1
\RequirePackage{ifpdf}
\ifpdf 
\documentclass[pdftex]{sigma}
\else
\documentclass{sigma}
\fi

\usepackage{rotating}

\usepackage{stmaryrd,mathtools}

\numberwithin{equation}{section}

\newtheorem{Theorem}{Theorem}[section]
\newtheorem*{Theorem*}{Theorem}
\newtheorem*{HypothesisI}{Hypothesis I}
\newtheorem*{HypothesisII}{Hypothesis II}

\newtheorem{Lemma}[Theorem]{Lemma}

\theoremstyle{definition}

\newtheorem{Remark}[Theorem]{Remark}

\newcommand{\artanh}{\operatorname{artanh}}
\newcommand{\even}{\operatorname{even}}
\newcommand{\odd}{\operatorname{odd}}
\newcommand{\GOE}{\operatorname{GOE}}
\newcommand{\OOE}{\operatorname{OE}}
\newcommand{\SE}{\operatorname{SE}}
\newcommand{\UE}{\operatorname{UE}}
\newcommand{\Ai}{\operatorname{Ai}}

\newcommand{\GUE}{\operatorname{GUE}}
\newcommand{\GSE}{\operatorname{GSE}}
\newcommand{\LUE}{\operatorname{LUE}}
\newcommand{\LSE}{\operatorname{LSE}}
\newcommand{\LOE}{\operatorname{LOE}}
\newcommand{\GbE}{\operatorname{G\beta E}}
\newcommand{\LbE}{\operatorname{L\beta E}}
\newcommand{\tr}{\operatorname{tr}}
\newcommand{\C}{{\mathbb C}}
\newcommand{\E}{{\mathbb E}}
\newcommand{\Q}{{\mathbb Q}}
\newcommand{\R}{{\mathbb R}}
\newcommand{\prob}{{\mathbb P}}

\begin{document}

\allowdisplaybreaks

\newcommand{\arXivNumber}{2506.18673}

\renewcommand{\PaperNumber}{047}

\FirstPageHeading

\ShortArticleName{Asymptotic Expansions at the Soft Edge III: Generating Functions}

\ArticleName{Asymptotic Expansions of Gaussian and Laguerre\\ Ensembles at the Soft Edge III: Generating Functions}

\Author{Folkmar BORNEMANN}

\AuthorNameForHeading{F.~Bornemann}

\Address{Department of Mathematics, Technical University of Munich, 80290 Munich, Germany}
\Email{\mail{bornemann@tum.de}}

\ArticleDates{Received January 08, 2026, in final form April 27, 2026; Published online May 13, 2026}

\Abstract{We conclude our work on asymptotic expansions at the soft edge for the classical $n$-dimensional Gaussian and Laguerre ensembles, now studying the gap-probability generating functions. We show that the correction terms in the asymptotic expansion are multilinear forms of the higher-order derivatives of the leading-order term, with certain rational polynomial coefficients that are {\em independent} of the dummy generating function variable. In this way, the same multilinear structure, with the same polynomial coefficients, is inherited by the asymptotic expansion of any linearly induced quantity such as the distribution of the $k$-th largest level. Whereas the results for the unitary ensembles are presented with proof, the discussion of the orthogonal and symplectic ones is based on some hypotheses. To substantiate the hypotheses, we check the result for the $k$-th largest level in the orthogonal ensembles against simulation data for choices of $n$ and $k$ that require as many as four correction terms to achieve satisfactory accuracy.}


\Classification{60B20; 15B52; 62E20; 41A60; 33E17}

\begin{flushright}
\begin{minipage}{63mm}
\em Dedicated to Peter Forrester\\ on the occasion of his 65th birthday
\end{minipage}
\end{flushright}

\section{Introduction}

\subsection{Generating functions and linearly induced quantities}
The unordered levels (eigenvalues) of the $n$-dimensional Gaussian and Laguerre random matrix ensembles,\footnote{We consider the orthogonal ($\beta=1$), unitary ($\beta=2$) and symplectic ($\beta=4$) ensembles -- parametrized as in Appendix~\ref{app:ensembles}.} denoted by $x_1,\dots,x_n$, constitute a finite simple point process on $\Omega\subset\R$ -- where $\Omega=(-\infty,\infty)$ in the Gaussian case and $\Omega=[0,\infty)$ in the Laguerre case. When focusing on the vicinity of the largest level, the soft edge, many important statistics are encoded by the generating function
\begin{equation*}
E(x;\xi) := \E\bigl((1-\xi)^{N(x,\infty)}\bigr), \qquad N(x,\infty) := \#\{j \mid x_j > x\}.
\end{equation*}
Note that $E(x;\xi)$ is a polynomial of degree $n$ in the dummy variable $\xi$.
Denoting the order statistic of the levels, from smallest to largest (there are almost surely no multiple levels), by
\[
x_{(1)} < x_{(2)} < \cdots < x_{(n)},
\]
we briefly list the following examples of quantities encoded by the generating function:\footnote{For details and the notion of `thinning' in (D), see Appendix~\ref{app:PGF}.}
\begin{itemize}
\item[(A)] the expected value of the number of levels larger than $x$ is
\[
\E(N(x,\infty)) = - \frac{\partial}{\partial\xi} E(x;\xi)\biggr|_{\xi=0},
\]
\item[(B)] the $k$-level gap probability is
\[
\prob (N(x,\infty)=k) = \frac{(-1)^k}{k!} \frac{\partial^k}{\partial\xi^k} E(x;\xi)\biggr|_{\xi=1},
\]
\item[(C)] the distribution of the $(k+1)$-th largest level $x_{(n-k)}$, $k=0,1,\dots,n-1$, is given by
\[
\prob\bigl(x_{(n-k)} \leq x\bigr) = \sum_{j=0}^k \prob(N(x,\infty)=j) = \sum_{j=0}^k \frac{(-1)^j}{j!} \frac{\partial^j}{\partial\xi^j} E(x;\xi)\biggr|_{\xi=1},
\]
\item[(D)] if $0\leq \xi\leq 1$, the value of the generating function has a probabilistic meaning in itself,
\[
E(x;\xi) = \prob (\text{$\max$-level of a $\xi$-thinning is $\leq x$} ).
\]
\end{itemize}
These examples share the feature that the quantity at hand, $E(x)$ say, is obtained by applying a~constant coefficient linear differential operator ${\mathcal L_\xi}\in \R[\partial_\xi]$ to $E(x;\xi)$, the result of which is then evaluated at some argument $\xi=\xi_*\in[0,1]$:
\begin{equation}\label{eq:E_induced}
E(x) = {\mathcal L_\xi} E(x;\xi) |_{\xi=\xi_*},
\end{equation}
we call such a quantity $E(x)$ {\em linearly induced} from the generating function $E(x;\xi)$.

We use indices as in $E_{\beta,n}$ to signify the dependence of a quantity $E$ on $\beta$ and $n$. In the Laguerre cases, we follow \cite[Section~2.2]{B25-3} and omit the {\em Wishart parameter} $p$ (defined in \eqref{eq:palpha}) from the notation: operations applied to $n$ are understood to apply in the same way to $p$. We recall that the Gaussian cases can be cast by (formally) taking $p\to \infty$ in the Laguerre ones.\footnote{
The reason is a $\LbE_{n,p}\to \GbE_n$ transition law as $p\to\infty$. See, e.g., \cite[Appendix B]{B25-1}, \cite{MR3413957} or, in the case of the $\LOE$ with integer $p$, \cite[Corollary 13.3.2]{MR1990662} and \cite[Corollary 9.5.7]{MR652932}.}

\subsection{Asymptotic expansions}
In generalizing our work \cite{B25-1,B25-3} on asymptotic expansions of the largest-level distributions and the level density at the soft edge as $n\to \infty$, we study in this paper expansions of the generating functions.

It is known that the largest level $x_{(n)}$ is asymptotically equal, almost surely, to a certain quantity $\mu_n$ when $n\to \infty$. The soft-edge scaling limit is based on an affine recentering and rescaling `zoom-in',
\[
x = \mu_n + \sigma_n s,
\]
with a certain rescaling parameter $\sigma_n$; the expansion is in powers of a parameter $h_n \asymp n^{-2/3}$~-- see equation~\eqref{eq:UEscaling} below for the precise expressions. Also, there is a parameter $0\leq \tau_n \leq 1$ satisfying $\tau_n=0$ in the Gaussian case and $\tau_n \to 0^+$ as~$p\to\infty$, for $n$ fixed, in the Laguerre case. To prevent the appearance of terms with odd powers of~$h^{1/2}$ when expanding in powers of the designated expansion parameter $h$, the parameters $\mu_n$, $4\sigma_n$, $h_n$, $\tau_n$ have to be modified for~${\beta=1,4}$ as per (which simultaneously applies to the `hidden' Wishart parameter~$p$ omitted from the notation)%
\begin{subequations}
\begin{equation}\label{eq:nprimeDef}
n' =
\begin{cases}
n - \frac12, & \beta = 1,\\
n, & \beta = 2,\\
2n+ \frac12, & \beta = 4,
\end{cases}
\end{equation}
noting the relation -- which mirrors the second relation in \eqref{eq:ForresterRains} below --
\begin{equation}\label{eq:nprime14}
(2n+1)'\big|_{\beta=1} = n'\big|_{\beta=4}.
\end{equation}
\end{subequations}
Our main result is the expansion
\begin{subequations}\label{eq:genTheory}
\begin{equation}\label{eq:genExpan}
E_{\beta,n}(x;\xi) |_{x=\mu_{n'} + \sigma_{n'} s} = F_\beta(s;\xi) + \sum_{j=1}^m G_{\beta,j}(s,\tau_{n'};\xi)h_{n'}^j + h_{n'}^{m+1}O({\rm e}^{-s}),
\end{equation}
which we obtain in Section~\ref{sect:unitary} for $\beta=2$ with proof (subject to $0\leq \xi \leq 1$ and other parameter constraints), and in Section~\ref{sect:orthsymp} for $\beta=1,4$ with a derivation based on certain structural hypotheses. Here, with certain polynomials $P_{\beta,j,k} \in \Q[s,\tau]$ that are {\em independent} of the dummy variable~$\xi$, the correction terms take the multilinear form
\begin{equation}
G_{\beta,j}(s,\tau;\xi) = \sum_{k=1}^{2j} P_{\beta,j,k}(s,\tau)\cdot \partial_s^k F_\beta(s;\xi).
\end{equation}
Moreover, there is a remarkable $\beta=1$ to $\beta=4$ duality that gives $P_{1,j,k}= P_{4,j,k}$.
\end{subequations}
Subject to the underlying constraints, the expansion \eqref{eq:genExpan} is uniformly valid for $0\leq \xi \leq 1$ as $s$ stays bounded from below; preserving uniformity, the expansion (including its remainder term) can be repeatedly differentiated with respect to the variables $s$ and $\xi$. This implies that any linearly induced quantity \eqref{eq:E_induced}, now written as
\begin{subequations}\label{eq:inducedTheory}
\begin{equation}
E_{n,\beta}(x) = {\mathcal L_\xi} E_{\beta,n}(x;\xi) |_{\xi=\xi_*},
\end{equation}
 enjoys the same form of expansion, namely
\begin{equation}
E_{\beta,n}(x) |_{x=\mu_{n'} + \sigma_{n'} s} = F_\beta(s) + \sum_{j=1}^m G_{\beta,j}(s,\tau_{n'})h_{n'}^j + h_{n'}^{m+1}O({\rm e}^{-s}).
\end{equation}
Here, by linearity, with exactly the {\em same} $\xi$-independent polynomials $P_{\beta,j,k}$ as in \eqref{eq:genTheory},
\begin{equation}\label{eq:GbetaMult}
G_{\beta,j}(s,\tau) = \sum_{k=1}^{2j} P_{\beta,j,k}(s,\tau)\cdot F_\beta^{(k)}(s),\qquad F_\beta(s) = {\mathcal L_\xi}F_\beta(s;\xi) |_{\xi=\xi_*}.
\end{equation}
\end{subequations}
That is, just the computational mastery of the leading term $F_\beta$ and its higher-order derivatives is required to evaluate {\em any} correction term in the asymptotic expansion of $E_{n,\beta}$. The coefficient polynomials $P_{\beta,j,k}$ are known: in \cite{B25-1}, when addressing the distribution of the largest level -- that is, the quantity $E_{\beta,n}(x;1)$ for which the limit laws $F_\beta(s)$ are the $\beta$-Tracy--Widom distributions -- we described an algorithm to compute the polynomials $P_{\beta,j,k}$. The concrete instances for up to $j=4$
are displayed\footnote{Their functional form quickly becomes bulky as the maximum total degree appears to grow like $3j$.} in Table~\ref{tab:1} for $\beta=2$, and in Table~\ref{tab:2} for $\beta=1,4$.

\begin{sidewaystable}
\centering
\vspace*{1mm}\small

\caption{The polynomial coefficients $P_{\beta,j,k}(s,\tau)$ for $\beta=2$ ($j=1,\dots,4$, $k=1,\dots,2j$); for the $\GUE$ take $\tau=0$.}\label{tab:1}
\vspace{1mm}\renewcommand{\arraystretch}{2.0}

$\begin{array}{@{}l@{}}
\hline
 P_{\beta,1,1}(s,\tau )=\dfrac{1}{5} s^2 (1-2 \tau ), \quad
 P_{\beta,1,2}(s,\tau )=\dfrac{1}{10}(\tau -3) \\[2mm]
 \hline
 P_{\beta,2,1}(s,\tau )=\dfrac{1}{175} s^3 \bigl(43 \tau ^2-18 \tau -8\bigr)+\dfrac{1}{350}
 \bigl(\tau ^2+94 \tau -141\bigr) \\[2mm]
 P_{\beta,2,2}(s,\tau )=\dfrac{1}{50} s^4 (2 \tau -1)^2+\dfrac{1}{175} s \bigl(-4 \tau ^2-26 \tau
 +39\bigr), \quad
 P_{\beta,2,3}(s,\tau )=-\dfrac{1}{50} s^2 (\tau -3) (2 \tau -1), \quad
 P_{\beta,2,4}(s,\tau )=\dfrac{1}{200} (\tau -3)^2 \\[2mm]\hline
P_{\beta,3,1}(s,\tau )=\dfrac{s^4 \bigl(-1384 \tau ^3+551 \tau ^2+212 \tau
 +148\bigr)}{7875}-\dfrac{2 s \bigl(4 \tau ^3-161 \tau ^2+1108 \tau -1108\bigr)}{7875} \\[2mm]
 P_{\beta,3,2}(s,\tau )=\dfrac{s^2 \bigl(12 \tau ^3-77 \tau ^2+328 \tau
 -265\bigr)}{1050}-\dfrac{s^5 (2 \tau -1) \bigl(43 \tau ^2-18 \tau -8\bigr)}{875} \\[2mm]
 P_{\beta,3,3}(s,\tau )=-\dfrac{s^6 (2 \tau -1)^3}{750} +\dfrac{s^3 \bigl(59 \tau ^3-51 \tau ^2-162
 \tau +102\bigr)}{1750}+\dfrac{61 \tau ^3+1351 \tau ^2-10403 \tau +10403}{31500} \\[2mm]
 P_{\beta,3,4}(s,\tau )=\dfrac{ s^4 (\tau -3) (2 \tau -1)^2}{500}-\dfrac{s (\tau -3) \bigl(4 \tau
 ^2+26 \tau -39\bigr)}{1750}, \quad
 P_{\beta,3,5}(s,\tau )=-\dfrac{s^2 (\tau -3)^2 (2 \tau -1)}{1000}, \quad
 P_{\beta,3,6}(s,\tau )=\dfrac{(\tau -3)^3}{6000} \\[2mm]\hline
 P_{\beta,4,1}(s,\tau )=\dfrac{2 s^5 \bigl(206303 \tau ^4-80156 \tau ^3-29708 \tau ^2-18264 \tau
 -14792\bigr)}{3031875}+\dfrac{s^2 \bigl(1751 \tau ^4-38732 \tau ^3+1834 \tau ^2+676152 \tau
 -720404\bigr)}{3031875} \\[2mm]
 P_{\beta,4,2}(s,\tau )=\dfrac{s^6 \bigl(55393 \tau ^4-48736 \tau ^3-1498 \tau ^2+1416 \tau
 +2648\bigr)}{551250}+\dfrac{s^3 \bigl(-9069 \tau ^4+38420 \tau ^3-35539 \tau ^2-253054 \tau
 +257512\bigr)}{1212750}\\
 \hphantom{P_{\beta,4,2}(s,\tau )=}{} +\dfrac{1519 \tau ^4-22028 \tau ^3+1809526 \tau ^2-6218452 \tau
 +4663839}{3465000} \\[2mm]
 P_{\beta,4,3}(s,\tau )=\dfrac{s^7 (2 \tau -1)^2 \bigl(43 \tau ^2-18 \tau -8\bigr)}{8750}+\dfrac{s^4
 \bigl(-6172 \tau ^4+7922 \tau ^3-7315 \tau ^2+22822 \tau -12843\bigr)}{220500}+\dfrac{s
 \bigl(-194 \tau ^4-332 \tau ^3-36251 \tau ^2+128972 \tau -96729\bigr)}{173250} \\[2mm]
 P_{\beta,4,4}(s,\tau )=\dfrac{s^8 (2 \tau -1)^4}{15000}-\dfrac{s^5 (2 \tau -1) \bigl(51 \tau ^3-99
 \tau ^2-58 \tau +63\bigr)}{8750}+\dfrac{s^2 \bigl(164 \tau ^4-3041 \tau ^3+32263 \tau
 ^2-61233 \tau +31405\bigr)}{220500} \\[2mm]
 P_{\beta,4,5}(s,\tau )=-\dfrac{s^6 (\tau -3) (2 \tau -1)^3}{7500}+\dfrac{s^3 (\tau -3) \bigl(15 \tau
 ^3+9 \tau ^2-74 \tau +36\bigr)}{7000}+\dfrac{(\tau -3) \bigl(113 \tau ^3+1883 \tau ^2-16999
 \tau +16999\bigr)}{630000} \\[2mm]
 P_{\beta,4,6}(s,\tau )=\dfrac{s^4 (\tau -3)^2 (2 \tau -1)^2}{10000}-\dfrac{s (\tau -3)^2 \bigl(4
 \tau ^2+26 \tau -39\bigr)}{35000}, \quad
 P_{\beta,4,7}(s,\tau )=-\dfrac{s^2 (\tau -3)^3 (2 \tau -1)}{30000}, \quad
 P_{\beta,4,8}(s,\tau )=\dfrac{(\tau -3)^4}{240000} \\[2mm]\hline
\end{array}$

\end{sidewaystable}

\begin{sidewaystable}
\centering
\vspace*{1mm}\small

\caption{The polynomial coefficients $P_{\beta,j,k}(s,\tau)$ for $\beta=1$ resp. $\beta=4$ ($j=1,\dots,4$, $k=1,\dots,2j$); for the $\GOE$ and $\GSE$ take $\tau=0$.}\label{tab:2}
\vspace{1mm}\renewcommand{\arraystretch}{2.0}

$\begin{array}{@{}l@{}}
\hline
P_{\beta,1,1}(s,\tau )=\dfrac{1}{5} s^2 (1-2 \tau ), \quad
 P_{\beta,1,2}(s,\tau )=\dfrac{1}{5}(\tau -3) \\[2mm]\hline
P_{\beta ,2,1}(s,\tau )=\dfrac{1}{175} s^3 \bigl(43 \tau ^2-18 \tau -8
\bigr)+\dfrac{1}{700} \bigl(9 \tau ^2+496 \tau -744\bigr) \\[2mm]
 P_{\beta ,2,2}(s,\tau )=\dfrac{1}{50} s^4 (2 \tau -1)^2-\dfrac{2}{175}
s \bigl(4 \tau ^2+26 \tau -39\bigr),\quad
 P_{\beta ,2,3}(s,\tau )=-\dfrac{1}{25} s^2 (\tau -3) (2 \tau -1),\quad
 P_{\beta ,2,4}(s,\tau )=\dfrac{1}{50} (\tau -3)^2 \\[2mm]\hline
 P_{\beta ,3,1}(s,\tau )=\dfrac{s^4 \bigl(-1384 \tau ^3+551 \tau ^2+212 \tau +148\bigr)}{7875}+\dfrac{s \bigl(-67 \tau ^3+1778 \tau ^2-12784 \tau +12784\bigr)}{15750} \\[2mm]
 P_{\beta ,3,2}(s,\tau )=\dfrac{s^2 \bigl(42 \tau ^3-449 \tau ^2+1600 \tau -1168\bigr)}{2100}-\dfrac{s^5 (2 \tau -1) \bigl(43 \tau ^2-18 \tau -8\bigr)}{875} \\[2mm]
 P_{\beta ,3,3}(s,\tau )=-\dfrac{ s^6 (2 \tau -1)^3}{750}+\dfrac{s^3 \bigl(59 \tau ^3-51 \tau ^2-162 \tau +102\bigr)}{875} +\dfrac{289 \tau ^3+6349 \tau ^2-46472 \tau +46472}{31500} \\[2mm]
 P_{\beta ,3,4}(s,\tau )=\dfrac{s^4 (\tau -3) (2 \tau -1)^2}{250} -\dfrac{2 s (\tau -3) \bigl(4 \tau ^2+26 \tau -39\bigr)}{875},\quad
 P_{\beta ,3,5}(s,\tau )=-\dfrac{s^2 (\tau -3)^2 (2 \tau -1)}{250} , \quad
 P_{\beta ,3,6}(s,\tau )=\dfrac{(\tau -3)^3}{750} \\[2mm]\hline
P_{\beta ,4,1}(s,\tau )=\dfrac{2 s^5 \bigl(206303 \tau ^4-80156 \tau ^3-29708 \tau ^2-18264 \tau -14792\bigr)}{3031875}+\dfrac{s^2 \bigl(27893 \tau ^4-427376 \tau ^3-224588 \tau ^2+8816736 \tau -9073472\bigr)}{12127500} \\[2mm]
 P_{\beta ,4,2}(s,\tau )=\dfrac{s^6 \bigl(55393 \tau ^4-48736 \tau ^3-1498 \tau ^2+1416 \tau +2648\bigr)}{551250}+\dfrac{s^3 \bigl(-29863 \tau ^4+222804 \tau ^3-82404 \tau ^2-1329280 \tau +1179296\bigr)}{2425500}\\
 \hphantom{P_{\beta ,4,2}(s,\tau )=}{}
 +\dfrac{349163 \tau ^4-2102656 \tau ^3+273755552 \tau ^2-945672704 \tau +709254528}{97020000} \\[2mm]
 P_{\beta ,4,3}(s,\tau )=\dfrac{s^7 (2 \tau \!-\!1)^2 \bigl(43 \tau ^2\!-\!18 \tau \!-\!8\bigr)}{8750}\!+\!\dfrac{s^4 \bigl(-24436 \tau ^4\!+\!37484 \tau ^3\!-\!44317 \tau ^2\!+\!101872 \tau \!-\!53640\bigr)}{441000}\!+\!\dfrac{s \bigl(-6367 \tau ^4\!-\!12211 \tau ^3\!-\!1141693 \tau ^2\!+\!4038016 \tau \!-\!3028512\bigr)}{1212750} \\[2mm]
 P_{\beta ,4,4}(s,\tau )=\dfrac{s^8 (2 \tau -1)^4}{15000}-\dfrac{s^5 (2 \tau -1) \bigl(51 \tau ^3-99 \tau ^2-58 \tau +63\bigr)}{4375}+\dfrac{s^2 \bigl(530 \tau ^4-14747 \tau ^3+143983 \tau ^2-265344 \tau +132424\bigr)}{220500} \\[2mm]
 P_{\beta ,4,5}(s,\tau )=-\dfrac{s^6 (\tau -3) (2 \tau -1)^3}{3750}+\dfrac{s^3 (\tau -3) \bigl(15 \tau ^3+9 \tau ^2-74 \tau +36\bigr)}{1750}+\dfrac{(\tau -3) \bigl(71 \tau ^3+1211 \tau ^2-10408 \tau +10408\bigr)}{45000} \\[2mm]
 P_{\beta ,4,6}(s,\tau )=\dfrac{s^4 (\tau -3)^2 (2 \tau -1)^2}{2500}-\dfrac{s (\tau -3)^2 \bigl(4 \tau ^2+26 \tau -39\bigr)}{4375},\quad
 P_{\beta ,4,7}(s,\tau )=-\dfrac{s^2 (\tau -3)^3 (2 \tau -1)}{3750}, \quad
 P_{\beta ,4,8}(s,\tau )=\dfrac{(\tau -3)^4}{15000}\\[2mm]\hline
\end{array}$

\end{sidewaystable}

\begin{figure}[t!]\centering
\includegraphics[width=0.485\textwidth]{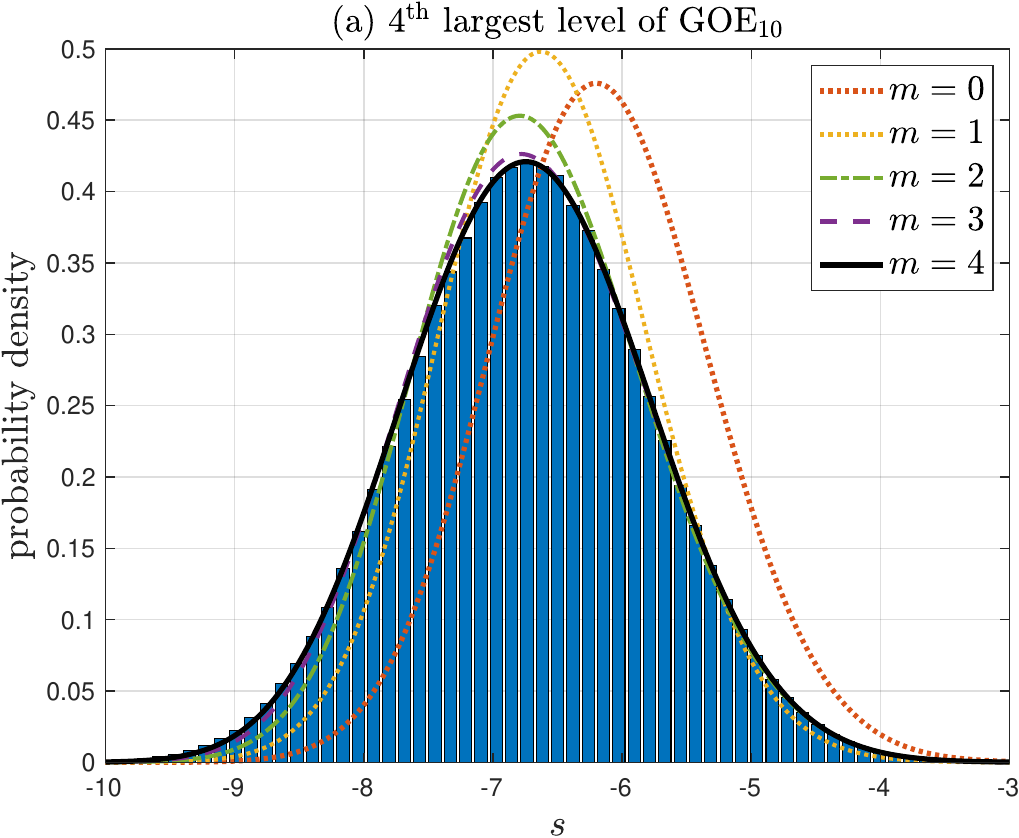}\hfil\,
\includegraphics[width=0.475\textwidth]{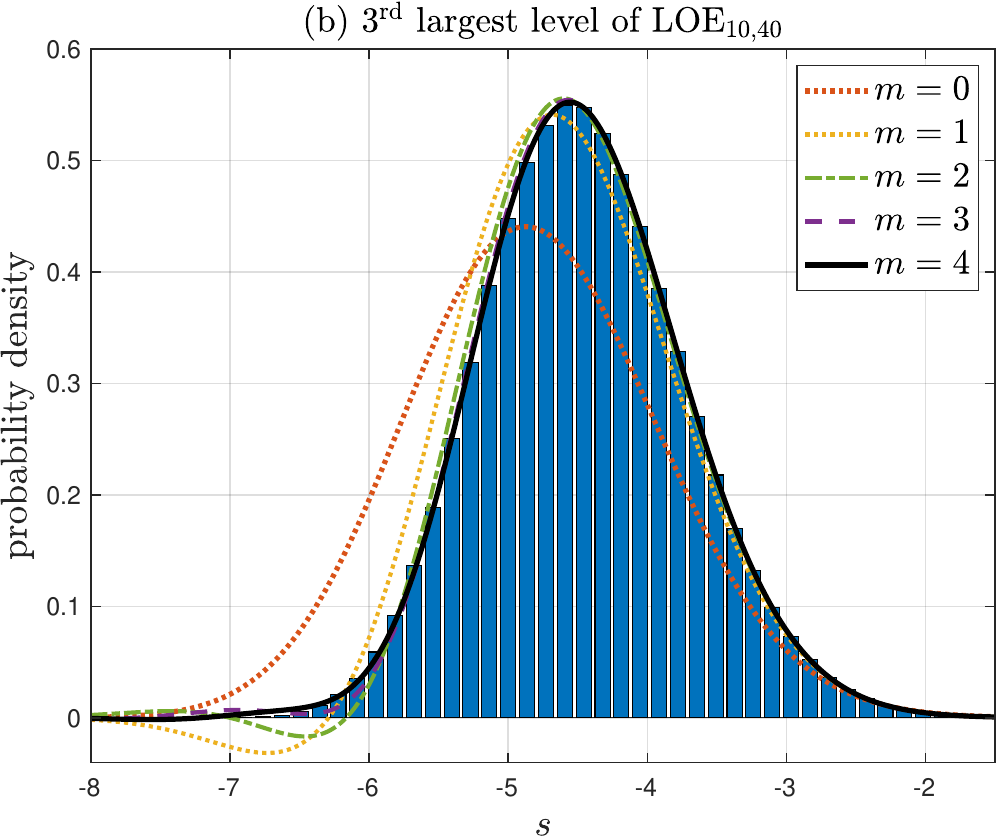}\hfil\,\\
\includegraphics[width=0.475\textwidth]{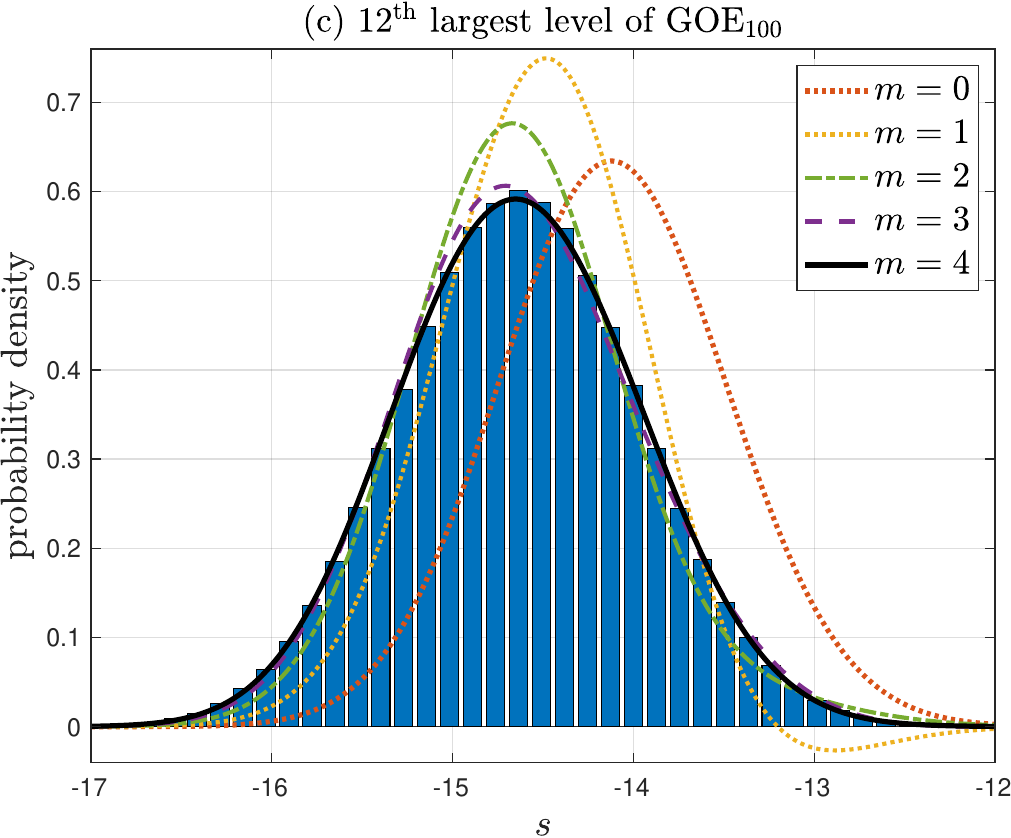}\hfil\,
\includegraphics[width=0.475\textwidth]{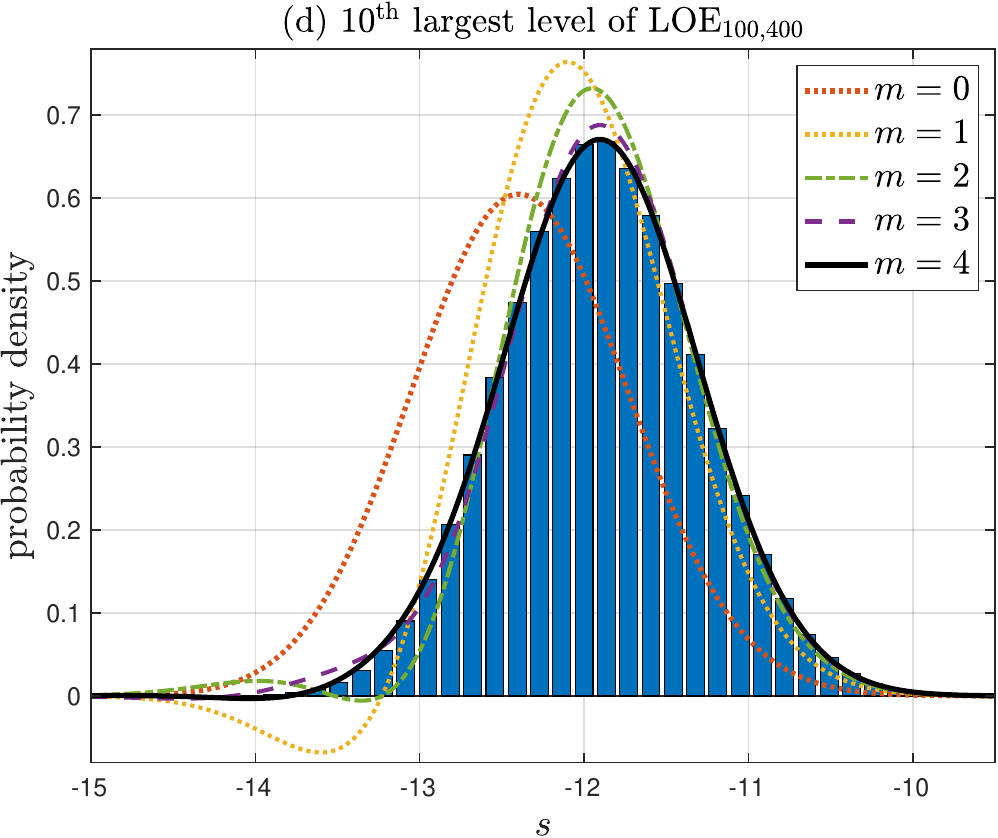}\hfil\,
\caption{{\footnotesize Histograms (blue bars) of the scaled $k^\text{th}$-largest level $(x_{(n-k+1)} - \mu_{n'})/\sigma_{n'}$ for simulations with $10^6$ draws from an orthogonal ensemble vs. the density approximations derived from the asymptotic expansion \eqref{eq:inducedTheory} with $m=0,1,2,3,4$. The indices $k$ were chosen to be relatively large compared to the dimension $n$ so that the limit law ($m=0$, red dotted line) is unsatisfactorily inaccurate, while it is not before $m=4$ (solid black line) that the approximation stays within the plotting accuracy. Note that intermediate approximations of the densities do not need to be positive. The limit laws and their derivatives were evaluated by the methods in \cite[Section~6]{MR2895091}; the polynomials $P_{1,j,k}$ were taken from~Table~\ref{tab:2}}.}\label{fig:evidence}
\end{figure}

\subsection{The evidence for the orthogonal and symplectic ensembles}\label{subsect:evidence} For $\beta=1,4$, the claimed asymptotic expansion \eqref{eq:genTheory} is based on some hypotheses and has thus the status of a conjecture only. Here is a list of some substantiating evidence and consistency.
\begin{itemize}\itemsep=0pt
\item The expansion \eqref{eq:genTheory} induces an expansion of the level densities
\[
\rho_{\beta,n}(x) = -\frac{\partial}{\partial x} \E(N(x,\infty)) = \frac{\partial^2}{\partial x\partial\xi} E(x;\xi)\biggr|_{\xi=0}.
\]
In our work \cite{B25-3},\footnote{Note added in proof: in the recent work \cite{2603.22974}, Forrester et al. discuss asymptotic expansions of certain scalar differential equations satisfied by the level densities at the soft and hard edge.} we establish asymptotic expansions of the level densities with {\em proof} even for $\beta=1,4$. The induced expansions agree with the directly proven ones for at least up to $m=10$. Moreover, using an algebraic independence result \cite{B25-2} for the Airy function, its derivative, and antiderivative, we could reconstruct from the proven expansion of the level densities the precise functional form of the presumed $12$ polynomials
\[
P_{\beta,1,1}, \qquad P_{\beta,1,2}, \qquad P_{\beta,2,1}, \qquad P_{\beta,2,2}, \qquad P_{\beta,2,3}, \qquad P_{\beta,2,4}, \qquad \beta=1,4,
\]
for orders $j=1$ and $j=2$;
see \cite[Section~5]{B25-3}.
\item The Forrester--Rains interrelations \eqref{eq:ForresterRains} readily imply the {\em interlacing property}
\[
\prob\bigl(x_{(n-k)} \leq x \mid \SE_n\bigr) = \prob\bigl(x_{(2n-2k)} \leq x \mid \OOE_{2n+1}\bigr),
\]
relating the order statistic of the orthogonal and symplectic ensemble. Since both distribution functions are linearly induced from the generating function, we get from~\eqref{eq:nprime14} that their respective asymptotic expansions must agree. That is, the limit laws are the same \cite[equation~(2.4)]{MR2181265} -- so comparing coefficients in \eqref{eq:GbetaMult}, assuming some linear independence, gives yet another derivation of the claimed duality $P_{1,j,k} = P_{4,j,k}$ which is independent of the one based on the explicit calculations in the proof of Theorem~\ref{thm:orthogonal} below.
\item Numerical experiments show the presumed asymptotic expansions to provide excellent approximations to simulation data already for comparatively small $n$. See Figure~\ref{fig:evidence} for the $4^\text{th}$-\big($12^\text{th}$-\big)largest level of $\GOE$ and the $3^\text{rd}$-($10^\text{th}$-)largest level of $\LOE$ with $p=4n$ for dimension $n=10$ ($n=100$) vs. the expansions evaluated up to $m=4$.\footnote{The distributions for finite $n$ can be evaluated numerically using the methods in \cite[Sections~2.2.3, 8.1 and~8.3]{MR2895091}, or \cite{2411.15635}.}
\item In \cite{2306.03798}, we addressed asymptotic expansions for the related hard-to-soft edge transition limit. Here, for $\xi=1$, we could establish the case $m=1$ with {\em proof}, showing that the corresponding coefficient polynomials for $\beta=1$ and $\beta=4$ agree. Further, up to $m=3$, we could confirm the structure of the expansion from highly accurate numerical experiments \cite[Appendix D.2]{2306.03798} -- including a precise rational reconstruction of the presumed polynomial coefficients, establishing their independence on $\xi$.
\end{itemize}

\subsection{Related work} For the soft-edge large matrix limits and the hard-to-soft edge transition limits, there has been a~considerable body of work on convergence rates and finite-size correction terms (i.e., the~case~${m=0}$ with a remainder term and the correction term for $m=1$). An extensive discussion of that literature can be found in our recent work \cite{arxiv:2301.02022,2306.03798,B25-1,B25-3}, where we establish the structure of the corresponding asymptotic expansions (with proof for $\beta=2$, subject to certain hypotheses for $\beta=1,4$). To our knowledge, there has been no prior study of convergence rates, finite-size corrections, or asymptotic expansions for the generating functions at the soft edge.

In the bulk, finite-size corrections of the generating function were studied for the classical circular ensembles in \cite{MR3647807}. There, the corresponding correction term for $m=1$ is given in operator-theoretic form and in terms of Painlev\'e transcendents. A much more telling form of those $m=1$ terms, namely as a polynomial coefficient (independent of the dummy variable) multiplying the second-order derivative of the limit law, was recently found by Forrester and Shen \cite[equations~(2.9), (3.6) and (3.12)]{ForrShen25}, who also consider asymptotic expansions -- however, they could not establish similarly compelling forms of the higher-order correction terms. We refer to their paper for more references on related work in the bulk and on the hard edge.

\section{Unitary ensembles}\label{sect:unitary}

The unitary Gaussian and Laguerre ensembles are known to be determinantal; cf., e.g., {\cite[Section~4.2]{MR2760897}}, \cite[Section~9]{MR2641363} and \cite[Section~19.1]{Mehta04}. In terms of the $L^2(\Omega)$-orthonormal Hermite and Laguerre wave functions $\phi_n$ -- induced by the orthogonal polynomials associated with the weights \eqref{eq:alphabeta} -- the correlation kernel for dimension $n$ is given by the finite rank projection kernel
\[
K_n(x,y) = \sum_{j=0}^{n-1} \phi_j(x)\phi_j(y).
\]
For brevity of notation in the Laguerre case, we omit here and in the following the dependence on the Wishart parameter $p$ as defined in~\eqref{eq:palpha}. The generating function can be expressed in terms of a Fredholm determinant (cf.~\eqref{eq:EAdet})
\begin{equation}\label{eq:E2n}
E_{2,n}(x;\xi) = \det(I-\xi K_n) |_{L^2(x,\infty)}.
\end{equation}
Now, the asymptotic expansion is obtained in complete analogy to the expansion of the distribution of the largest level \cite[Theorems~2.1, 3.1]{B25-1}, which amounts to taking $\xi=1$.
The scaling and expansion parameters introduced in \cite[Theorems~2.1, 3.1]{B25-1} are as follows:
\begin{subequations}\label{eq:UEscaling}
\begin{itemize}\itemsep=0pt
\item $\GUE_n$:
\begin{equation}\label{eq:GUEscaling}
\mu_n=\sqrt{2n}, \qquad \sigma_n = 2^{-1/2}n^{-1/6}, \qquad h_n =\frac{\sigma_n}{2\mu_n}= 4^{-1} n^{-2/3}, \qquad \tau_n=0.
\end{equation}
Note that $h_n \asymp n^{-2/3}$ as $n\to\infty$.
\item $\LUE_{n,p}$:
\begin{alignat}{3}
&\mu_{n} = \bigl(\sqrt{n}+\sqrt{p}\bigr)^2,\qquad && \sigma_{n}= \bigl(\sqrt{n}+\sqrt{p}\bigr) \bigg(\frac{1}{\sqrt{n}} + \frac{1}{\sqrt{p}}\bigg)^{1/3},&\nonumber\\
&h_{n} = \frac{\sigma_{n}}{\tau_{n}\,\mu_{n}}= \frac{1}{4} \bigg(\frac{1}{\sqrt{n}} + \frac{1}{\sqrt{p}}\bigg)^{4/3},\qquad&&
\tau_{n} = \frac{4}{\bigl(\sqrt{n}+\sqrt{p}\bigr) \bigl(\frac{1}{\sqrt{n}} + \frac{1}{\sqrt{p}}\bigr)}.&\label{eq:LUEscaling}
\end{alignat}
Note that $0 < \tau_n \leq 1$ and $p > n-1$, so that $p\to\infty$ and $h_n \asymp n^{-2/3}$ as $n\to\infty$.
\end{itemize}
\end{subequations}
With these parameters at hand, the results for the Gaussian and the Laguerre case can be stated as a single theorem.

\begin{Theorem}\label{thm:unitary} There holds, for any fixed non-negative integer $m \leq m_*$,\footnote{\label{fn:mstar}We restrict to $m\leq m_*$ for two {\em independent} reasons: Step 1 of the proof relies on a verification that the kernel expansion yields certain finite-rank correction kernels, while Step 3 hinges on the solvability of certain over\-determined linear systems over $\Q[s]$; both issues can currently be settled only by an explicit casewise inspection (however, a Riemann--Hilbert analysis, as in \cite{arXiv:2309.06733}, would probably dispense with the first one). Using a computer algebra system, we have confirmed the necessary conditions for all $m\leq 10$.}
\begin{subequations}
\begin{equation}\label{eq:E2expand}
E_{2,n}(x;\xi) |_{x=\mu_{n} + \sigma_{n} s} = F_2(s;\xi) + \sum_{j=1}^m G_{2,j}(s,\tau_{n};\xi)h_{n}^j + h_{n}^{m+1}O\bigl({\rm e}^{-2s}\bigr),\qquad n\to\infty,
\end{equation}
uniformly for $0\leq \xi \leq 1$ while $s$ stays bounded from below and, in the Laguerre case, $\tau_n$~stays bounded away from zero -- the Gaussian case being encoded with $\tau_n=0$. Preserving uniformity, the expansion can be repeatedly differentiated with respect to the variables $s$ and $\xi$.
The leading term is given by the Airy kernel determinant\/\footnote{For $n\to\infty$, and a fixed Laguerre parameter $\alpha = p-n>-1$ in the Laguerre case, the pointwise limit
\[
\lim_{n\to \infty} E_{2,n}(\mu_{n} + \sigma_{n} s;\xi) = \det(I- \xi K_{\Ai}) |_{L^2(s,\infty)}
\]
was established by Forrester \cite{MR1236195}; see also \cite[Sections~7.1, 7.2 and 9.1]{MR2641363}. $\Ai$ is the standard Airy function \cite[equation~(9.5.1)]{MR2723248}.}
\begin{equation}
F_2(s;\xi) = \det(I- \xi K_{\Ai}) |_{L^2(s,\infty)},\qquad K_{\Ai}(x,y) = \frac{\Ai(x)\Ai'(y)-\Ai'(x)\Ai(y)}{x-y},
\end{equation}
and the correction terms take the multilinear form
\begin{equation}\label{eq:G2struct}
G_{2,j}(s,\tau;\xi) = \sum_{k=1}^{2j} P_{2,j,k}(s,\tau)\cdot \partial_s^k F_2(s;\xi)
\end{equation}
\end{subequations}
with certain polynomials $P_{2,j,k} \in \Q[s,\tau]$ that are independent of the dummy variable~$\xi$ -- the concrete instances up to $j=4$ are displayed in Table~{\rm \ref{tab:1}}.
\end{Theorem}
\begin{proof} The proof is basically the same as the one given for $\xi=1$ in \cite[Theorems~2.1 and~3.1]{B25-1}, stated there separately for the Gaussian and the Laguerre cases. While recalling the principal steps here, we focus on the dependence on the dummy variable $\xi$.\footnote{The uniformity of the expansion for
$0\leq \xi \leq 1$ and its repeated differentiability with respect to $\xi$ follow from analyticity properties, which allow us to apply Ritt's theorem \cite[Theorem~1.4.2]{MR0435697}: $E_{2,n}(x;\xi)$ is a polynomial in~$\xi$ and the limit generating function $F_2(s;\xi)$, as a Fredholm determinant, is entire with respect to $\xi$.}

{\em Step 1: Kernel expansion.}
In \cite[Lemmas~2.1, 2.2, 3.1 and 3.2]{B25-1}, we proved from the asymptotic expansions of the Hermite and Laguerre wave functions \cite[Section~10]{B25-1} that, for $0\leq m \leq m_*$,
\begin{align*}
&\sigma_n K_n(\mu_n+\sigma_n x,\mu_n+\sigma_n y) \nonumber\\
&\qquad =  K_{\Ai}(x,y) +  \sum_{j=1}^m L_j(x,y;\tau_n) h_n^j + h_n^{m+1}O\bigl({\rm e}^{-(x+y)}\bigr),
\qquad n\to\infty,
\end{align*}
uniformly for $x$, $y$ bounded from below and, in the Laguerre case, as $\tau_n$ stays bounded away from zero. Here, the correction kernels $L_j$ are symmetric finite rank kernels of the form
\begin{equation*}
L_j(x,y;\tau) = \sum_{\mu=0,\nu=0}^{\lambda_j} {\mathfrak p}_{j,\mu\nu}(\tau) \Ai^{(\mu)}(x) \Ai^{(\nu)}(y)
\end{equation*}
with certain integers $\lambda_j$ and coefficient polynomials ${\mathfrak p}_{j,\mu\nu} \in \Q[\tau]$.

{\em Step 2: Lifting to the Fredholm determinant.} By the perturbation theory of determinants (see~\cite[Section~2.2]{arxiv:2301.02022}, \cite[Appendix~A]{2306.03798} and \cite[Appendix~A]{B25-1}), the kernel expansion inserted to~\eqref{eq:E2n} lifts to
\[
E_{2,n}(x;\xi) |_{x=\mu_{n} + \sigma_{n} s} = F_2(s;\xi) + F_2(s;\xi) \sum_{j=1}^m d_j(s,\tau_n;\xi) h_n^j + h_n^{m+1} O\bigl({\rm e}^{-2s}\bigr)
\]
with
\begin{equation}\label{eq:djT}
d_j(s,\tau;\xi) = \sum_{T} p_{j,T}(\tau)\cdot T(s;\xi).
\end{equation}
Here, the sum is taken over all terms $T=T(s;\xi)$ of the form
\begin{equation}\label{eq:T_term}
T = \begin{vmatrix}
u_{\mu_1\nu_1} & \cdots & u_{\mu_1\nu_\kappa} \\
\vdots & & \vdots \\
u_{\mu_\kappa\nu_1} & \cdots & u_{\mu_\kappa\nu_\kappa}
\end{vmatrix}, \qquad u_{\mu\nu} = \xi \tr\bigl((I-\xi K_{\Ai})^{-1} \Ai^{(\mu)}\otimes \Ai^{(\nu)}\bigr)\big|_{L^2(s,\infty)},
\end{equation}
with certain coefficient polynomials $p_{j,T}\in \Q[\tau]$ -- of which, for a given $j$, just finitely many are non-zero. That is, we obtain the expansion \eqref{eq:E2expand} with
\[
G_{2,j}(s,\tau;\xi) = F_2(s;\xi) d_j(s,\tau;\xi),
\]
and it remains to convert those highly nonlinear operator-theoretic expressions\footnote{Which, however, would be amenable for numerical evaluation by the methods detailed in \cite{MR2895091} and~\cite[Appendix]{MR3647807}.} to the much more telling multilinear form \eqref{eq:G2struct}.

{\em Step 3: Applying Tracy--Widom theory.} The quantities $u_{\mu\nu}$ and the limit generating function
\[
F_2(s;\xi) = \det(I- \xi K_{\Ai}) |_{L^2(s,\infty)}
\]
can be obtained as follows: replace in the defining expressions for $\xi=1$ all appearances of the Airy function $\Ai$ by $\sqrt{\xi}\Ai$. The same structure applies to the Tracy--Widom theory~\cite{MR1257246}, which is based on purely algebraic manipulations of operator determinants by using just the Airy differential equation -- while guaranteeing the existence of terms and expressions from the superexponential decay of $\Ai(s)\to 0$ as $s\to \infty$. In particular, \cite[equations~(1.11)--(1.15)]{MR1257246}
 shows that the logarithmic derivative of $F_2(s;\xi)$ with respect to $s$ is represented in terms of a~Painlev\'e~II transcendent $q(s;\xi)$: writing \smash{\raisebox{-0.5pt}{$F_2^{(k)}=\partial_s^k F_2(s;\xi)$}} and \smash{$ q^{(k)}=\partial_s^k q(s;\xi)$} for brevity, we~have
\begin{gather}
u_{00} = F_2'/F_2 = q'^2 - s q^2 - q^4, \nonumber\\
q'' = s q + 2 q^3,\qquad q = q(s;\xi) \sim \sqrt{\xi} \Ai(s),\qquad s\to\infty.\label{eq:PII}
\end{gather}
That is, $q$ is, for $\xi=1$, the Hastings--McLeod \cite{MR555581} solution of Painlev\'e~II and, for $0\leq \xi<1$, the Ablowitz--Segur~\cite{zbMATH05784617} solution. Note that $\xi$ enters \eqref{eq:PII} only via the boundary condition of $q$ as~${s\to\infty}$. Since the assertion of the theorem is trivial if $\xi=0$ (because then $E_{2,n}\equiv F_2 \equiv 1$), we can assume $0< \xi \leq 1$ from now on, which implies $q\not\equiv 0$.

The same underlying algebraic structure applies to the extension of the Tracy--Widom theory developed by Shinault and Tracy \cite{MR2787973} to study the quantities $u_{\mu\nu}$ -- a fortiori, it applies also to the algorithmic framing in \cite[Appendix~B]{arxiv:2301.02022}. The algorithm is based on the following observation: by Painlev\'e~II, the ring $\Q[s,q,q']$ is closed under differentiation with respect to the variable $s$ and we get inductively, by repeated differentiation of \eqref{eq:PII}, that
\begin{equation}\label{eq:F2higher}
F_2^{(k)}/F_2 \in \Q\bigl[s,q,q'\bigr];
\end{equation}
for instance
\begin{equation}\label{eq:F22nd}
F_2''/F_2 = -2 q^4 q'^2-2 s q^2 q'^2+q'^4+s^2 q^4+q^8+2 s q^6-q^2.
\end{equation}
Here, as for $k=1$ in \eqref{eq:PII}, the polynomial representation in $\Q[s,q,q']$ is {\em independent} of~$\xi$: the functional dependence on $\xi$ enters only indirectly via the boundary condition of the func\-tion~$q$ as $s\to\infty$. Since it is known for Painlev\'e~II (cf., e.g., \cite[Theorem~21.1]{MR1960811}) that $0\not\equiv q,q'$ are algebraically independent over $\C(s)$, we can handle $s$, $q$, $q'$ as three indeterminates when working with $\Q[s,q,q']$: equality in $\Q[s,q,q']$ is equivalent to equality of coefficients.

For each term $T$ of the form \eqref{eq:T_term}, the algorithm of \cite[Appendix~B]{arxiv:2301.02022} starts with the representations \eqref{eq:F2higher} and, whenever certain overdetermined linear systems over the ring $\Q[s]$ are solvable, computes representations $u_{\mu\nu}\in \Q[s,q,q']$ and, a fortiori, $T\in \Q[s,q,q']$. In a follow-up step, the algorithm returns the polynomial coefficients $r_{1,T},\dots,r_{\omega_T,T} \in \Q[s]$ appearing in the multilinear ansatz\footnote{Actually, it turns out that $\omega_T = \mu_1+\cdots +\mu_\kappa + \nu_1 + \cdots +\nu_\kappa + \kappa$, which we call the {\em order} of the term $T$.}
\begin{equation}\label{eq:ansatz}
T = r_{1,T} \frac{F_2'}{F_2} + \cdots + r_{\omega_T,T}\frac{F_2^{(\omega_T)}}{F_2}
\end{equation}
if satisfiable. Namely, by comparing the coefficients of monomials in $q$, $q'$ in that ansatz, which is an equation in the polynomial ring $\Q[s,q,q']=\Q[s][q,q']$, it can be
rendered as an overdetermined linear system over $\Q[s]$. By construction, that linear system, the question of its solvability, and the presumed polynomial solutions~$(r_{\lambda,T})_\lambda$ are all {\em independent} of $\xi$. In fact, when $j\leq m_*$, all the relevant linear systems are {\em uniquely} solvable.\footnote{Uniqueness in general would be equivalent to the sequence $F_2'/F_2,F''_2/F_2,F'''_2/F_2,\dots$ being linearly independent over $\Q[s]$.} While we conjecture that this holds in general, it remains an open problem whether the question can be addressed by any method other than checking, case by case, the (unique) solvability of the linear systems using symbolic software.

{\em Step 4: Wrap-up.} For at least $j\leq m_*$, by inserting the representations \eqref{eq:ansatz} of the terms $T$ into \eqref{eq:djT}, we get
\[
\begin{gathered}
G_{2,j}(s,\tau;\xi) = \sum_{k=1}^{\varpi_j} P_{2,j,k}(s,\tau)\cdot \partial_s^k F_2(s;\xi), \qquad P_{2,j,k}(s,\tau) = \sum_T r_{k,T}(s) p_{j,T}(\tau).
\end{gathered}
\]
Actually, for $j\leq m_*$, all the terms in the last sum cancel for orders $k > 2j$ so that we effectively get $\varpi_j = 2j$ \big(which is, for larger $j$, considerably smaller than the maximum order of some of the intermediate terms $T$ which enter the calculation\footnote{This is the reason why we took the detour of using certain nonlinear simplifying transformations in \cite[Sections~2.2 and 3.2]{B25-1} to shortcut the computation of the polynomials $P_{2,j,k}$.}\big).
\end{proof}

\begin{Remark} We expect that the bound $m_*$ in Theorem~\ref{thm:unitary} is insignificant and can be removed. Since the Gaussian case corresponds to taking $\tau_n=0$ and there is the $\LUE_{n,p}\to \GUE_n$ transition law for fixed $n$ and $p\to \infty$ (which implies $\tau_n\to 0$), we conjecture that the Laguerre case of Theorem~\ref{thm:unitary} holds uniformly for all $0< \tau_n\leq 1$: this would allow us to infer the Gaussian case of the theorem from the Laguerre one by simply exchanging the order of limits; cf.\ also \cite[Remarks~3.1 and 10.5]{B25-1} and \cite[Remark~4.1]{B25-3}.
\end{Remark}

\section{Orthogonal and symplectic ensembles}\label{sect:orthsymp}

\subsection{Interrelations of the three symmetry classes}
As in \cite[Section~4]{B25-1} for the largest level, the expansion terms of the orthogonal and symplectic ensembles will be obtained by observing a tight algebraic relation to the expansion terms of the unitary ensemble. The starting point are the interrelations, through decimation and superposition, between the three symmetry classes as established by Forrester and Rains \cite[Sections~4 and~5]{MR1842786} (see also \cite[Theorem~2]{MR2461989}) -- namely
\begin{alignat*}{3}
&\GUE_n = \even(\GOE_{n} \cup \GOE_{n+1}), \qquad&& \LUE_n = \even(\LOE_{n,p} \cup \LOE_{n+1,p+1}),& \\
&\GSE_n = \even(\GOE_{2n+1}),\qquad&& \LSE_{n,p} = \even(\LOE_{2n+1,2p+1}).&
\end{alignat*}
These relations state that the ordered levels of the ensembles on the left are distributed as every second level (`decimation') of the corresponding cases on the right, the union (`superposition') being taken over stochastically independent orthogonal ensembles.

We omit the explicit dependence on $p$ from the notation. It will tacitly be understood that whenever the dimension $n$ gets modified, the Wishart parameter $p$ is modified in {\em exactly the same way} (called the {\em index consistency rule} in \cite[Section~2.2]{B25-3}). That being said, we briefly~write
\begin{equation}\label{eq:ForresterRains}
\UE_n = \even(\OOE_{n} \cup \OOE_{n+1}), \qquad \SE_n = \even(\OOE_{2n+1}).
\end{equation}
We aim to restate these interrelations in terms of generating functions. In a first step, briefly writing the $k$-level gap probabilities as
\[
{\mathfrak E}_{\beta,n}(k; x) := \prob (N_{\beta,n}(x,\infty) = k),
\]
we infer from \eqref{eq:ForresterRains} the recursions
\begin{gather}
{\mathfrak E}_{2,n}(k; x) = \sum_{j=0}^{2k} {\mathfrak E}_{1,n}(j; x){\mathfrak E}_{1,n+1}(2k-j; x) + \sum_{j=0}^{2k+1} {\mathfrak E}_{1,n}(j; x){\mathfrak E}_{1,n+1}(2k+1-j; x),\nonumber\\
{\mathfrak E}_{4,n}(k; x) = {\mathfrak E}_{1,2n+1}(2k;x) + {\mathfrak E}_{1,2n+1}(2k+1;x),\label{eq:Erec}
\end{gather}
cf.\ the arguments leading to \cite[equation~(6.16)]{MR2895091}. We introduce the generating functions
\begin{equation*}
E_{\left\{\even\atop\odd\right\},n}(x;\xi) :=\sum_{k=0}^\infty (1-\xi)^{k} {\mathfrak E}_{1,n}\Bigl({\textstyle\left\{2k\atop2k+1\right\}}; x\Bigr)
\end{equation*}
and recall that
\[
E_{\beta,n}(x;\xi) =\sum_{k=0}^\infty (1-\xi)^k {\mathfrak E}_{\beta,n}(k; x);
\]
all of which are polynomials in $\xi$. Now, the recursions \eqref{eq:Erec} can be re-expressed as
\begin{subequations}\label{eq:EgenRel}
\begin{gather}
E_{2,n}(x;\xi) = E_{\even,n}(x;\xi)E_{\even,n+1}(x;\xi) + E_{\even,n}(x;\xi)E_{\odd,n+1}(x;\xi)\nonumber\\
 \hphantom{E_{2,n}(x;\xi) =}{}
 + E_{\odd,n}(x;\xi)E_{\even,n+1}(x;\xi) + (1-\xi) E_{\odd,n}(x;\xi)E_{\odd,n+1}(x;\xi),\nonumber\\
E_{4,n}(x;\xi) = E_{\even,2n+1}(x;\xi) + E_{\odd,2n+1}(x;\xi).\label{eq:EgenRel24}
\end{gather}
On the other hand, if we define a transformed dummy variable $\bar\xi$ by
\[
1-\bar\xi = (1-\xi)^2, \qquad \bar\xi := \xi(2-\xi),
\]
we have, by construction,
\begin{equation}\label{eq:EgenRel1}
E_{1,n}(x;\xi)= E_{\even,n}\bigl(x;\bar \xi\,\bigr) + (1-\xi) E_{\odd,n}\bigl(x;\bar\xi\,\bigr).
\end{equation}
\end{subequations}

These interrelations of the generating functions can be recast in a form that is tailored for generalizing the algebraic method of \cite[Section~4.2]{B25-1} -- see footnote~\ref{fn:algStruc}:

\begin{Lemma} Consider $0\leq \xi \leq 1$ and define the
functions\footnote{The definition of $E_{\pm,n}(x;\xi)$ mirrors an analogous construction \cite[equation~(6.24)]{MR2895091} in the soft-edge scaling limit.}
\begin{equation}\label{eq:Epm}
E_{\pm,n}(x;\xi) := E_{\even,n}(x;\xi) + \bigl(1\mp\xi^{1/2}\bigr) E_{\odd,n}(x;\xi),
\end{equation}
which are polynomials in $\xi^{1/2}$.
Then, the generating functions $E_{\beta,n}$ can be expressed in terms of~$E_{\pm,n}$ as follows:
\begin{subequations}\label{eq:Ebeta_Epm}
\begin{align}
&E_{1,n}(x;\xi)= \frac{1}{2} E_{+,n}\bigl(x;\bar\xi\,\bigr) \bigl(1+\xi/{\bar\xi}^{\,1/2}\bigr)+ \frac{1}{2} E_{-,n}\bigl(x;\bar\xi\,\bigr) \bigl(1-\xi/{\bar\xi}^{\,1/2}\bigr),\label{eq:E1Epm}\\
&E_{2,n}(x;\xi)= \frac{1}{2} E_{+,n}(x;\xi) E_{-,n+1}(x;\xi) + \frac{1}{2} E_{+,n+1}(x;\xi) E_{-,n}(x;\xi),\label{eq:E2Epm}\\
&E_{4,n}(x;\xi)= \frac{1}{2} E_{+,2n+1}(x;\xi) + \frac{1}{2} E_{-,2n+1}(x;\xi).\label{eq:E4Epm}
\end{align}
\end{subequations}
Here, observing \smash{$\xi/\bar\xi^{\,1/2}= \sqrt{\xi/(2-\xi)}$}, the right-hand sides are well defined for all $0\leq \xi\leq 1$.
\end{Lemma}
\begin{proof} A routine calculation shows that inserting \eqref{eq:Epm} into \eqref{eq:Ebeta_Epm} reclaims \eqref{eq:EgenRel}.
\end{proof}

\subsection{Soft-edge scaling limit} Using the same rescaling parameters $\mu_n$, $\sigma_n$ as in the $\beta=2$ case -- displayed in \eqref{eq:UEscaling} -- it is known that pointwise\footnote{See, e.g., \cite{MR2895091,MR2181265,MR2275509}; in the Laguerre case, the results are given in the literature for a fixed Laguerre parameter $\alpha = p-n>-1$.}
\begin{equation}\label{eq:SESL}
\lim_{n\to\infty} E_{\beta,n}(x;\xi) |_{x=\mu_{n_*}+\sigma_{n_*} s} = F_\beta(s;\xi), \qquad 0\leq \xi\leq 1 ,
\end{equation}
where $n_*=n$ for $\beta=1,2$ and $n_* = 2n$ for $\beta=4$. Using the differentiability of $F_\beta(s;\xi)$, we see that the same remains true if we replace
$n_*$ by any constant shift $n_* +c$.

For $0\leq \xi\leq 1$, Dieng \cite[equations~(2.1) and (2.2)]{MR2181265} and Forrester \cite[Corollary~1]{MR2275509} have found the following expressions for the limit functions in terms of Painlev\'e~II and Fredholm determinants -- see also \cite[equations~(6.5), (6.10) and (6.27)]{MR2895091}:
\begin{subequations}\label{eq:Fbeta_Fpm}
\begin{gather}
F_1(s;\xi)= \frac{1}{2} F_{+}\bigl(s;\bar\xi\,\bigr) \bigl(1+\xi/{\bar\xi}^{\,1/2}\bigr)+ \frac{1}{2} F_{-}\bigl(s;\bar\xi\,\bigr) \bigl(1-\xi/{\bar\xi}^{\,1/2}\bigr),\label{eq:F1Fpm}\\
F_2(s;\xi)= F_{+}(s;\xi) F_{-}(s;\xi) ,\label{eq:F2Fpm}\\
F_4(s;\xi)= \frac{1}{2} F_{+}(s;\xi) + \frac{1}{2} F_{-}(x;\xi),\label{eq:F4Fpm}\\
\intertext{with the auxiliary functions}
 F_{\pm}(s;\xi) = \sqrt{F_2(s;\xi)} \exp\biggl(\mp \frac{1}{2}\int_s^\infty q(t;\xi)\,{\rm d}t \biggr) = \det\bigl(I\mp\xi^{1/2} V_{\Ai}\bigr)\big|_{L^2(s,\infty)}.\label{eq:Fpm}
\end{gather}
\end{subequations}
Here, $q(s;\xi)$ denotes the Painlev\'e~II transcendent defined in \eqref{eq:PII}, and $V_{\Ai}$ is the integral operator with kernel
\[
V_{\Ai}(x,y) = \frac{1}{2}\Ai\biggl(\frac{x+y}{2}\biggr).
\]
In particular, the representation as a Fredholm determinant shows that $\xi \mapsto F_{\pm}(s;\xi)$ behaves analytically like an entire function applied to $\xi^{1/2}$.

Comparing \eqref{eq:Ebeta_Epm} with \eqref{eq:Fbeta_Fpm} strongly suggests the following lemma. Its proof requires solving both sets of equations for $E_\pm(x;\xi)$ and $F_\pm(s;\xi)$, though.

\begin{figure}[t!]\centering
\includegraphics[width=0.5\textwidth]{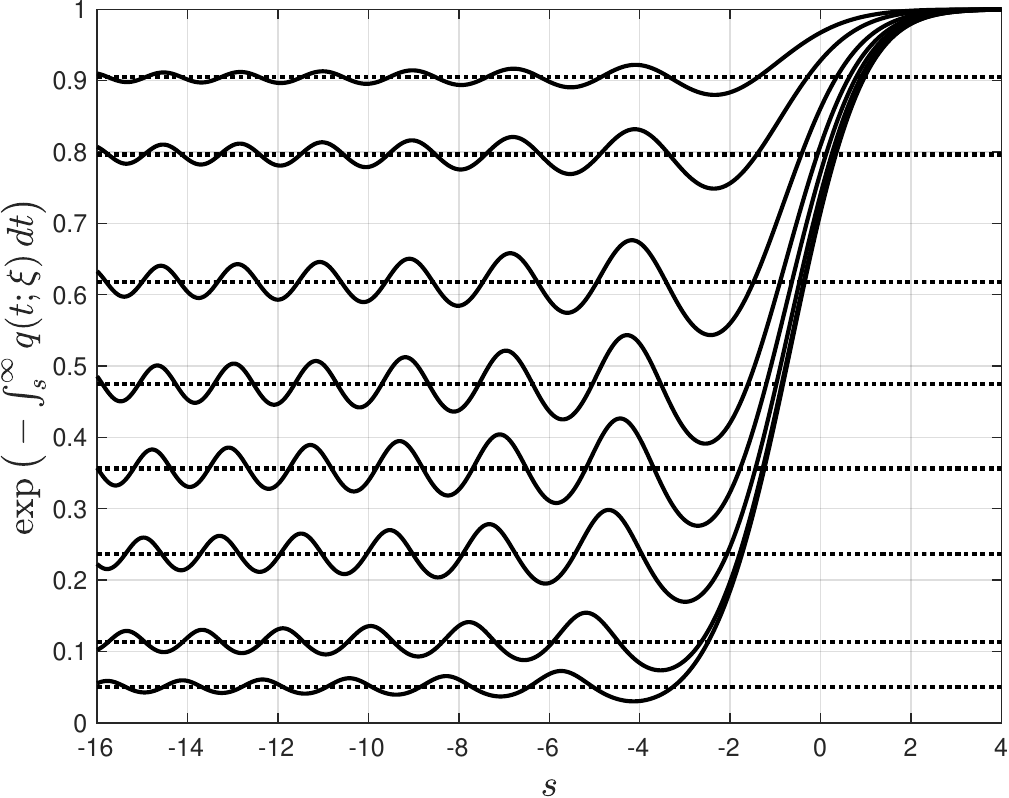}
\caption{{\footnotesize Plot of $\exp\bigl(-\int_{s}^\infty q(t;\xi)\,{\rm d}t\bigr)$ (solid black lines) vs. $\exp\bigl(-\artanh\xi^{1/2}\bigr)$ (dotted black lines) for $\xi = 0.01$ (topmost), $0.05, 0.2, 0.4, 0.6, 0.8, 0.95, 0.99$ (bottommost). Note that, as stated in \eqref{eq:total}, the solid lines converge to the dotted ones as $s\to-\infty$.}}
\label{fig:2}
\end{figure}

\begin{Lemma}\label{lem:EpmFpm} There is
\begin{equation}\label{eq:EpmFpm}
\lim_{n\to\infty} E_{\pm,n}(x;\xi) |_{x=\mu_{n}+\sigma_{n} s} = F_\pm(s;\xi).
\end{equation}
If $n$ is even and $0<\xi<1$, the proof requires the exclusion\footnote{If the limits are uniform, the singularities causing the exclusion points are removable.} of a discrete set $($which depends on~$\xi)$ of values of $s$.
\end{Lemma}
\begin{proof}
The assertion is clear for $\xi=0$ since then
\[
E_{\pm,n}(x;0) = E_{\even,n}(x;0) + E_{\odd,n}(x;0) = E_{1,n} (x;0) \equiv 1, \qquad F_{\pm}(s;0) \equiv 1.
\]
We now assume $0<\xi\leq 1$. Let $\xi_*:=1-\sqrt{1-\xi}$, so that $\bar\xi_* = \xi$. The equations \eqref{eq:E1Epm}/\eqref{eq:E4Epm} give the linear system
\[
\begin{pmatrix}
E_{1,2n+1}(x,\xi_*)\\
E_{4,n}(x,\xi)
\end{pmatrix}
=
\frac{1}{2}\begin{pmatrix}
1+\xi_*/\xi^{1/2} & 1-\xi_*/\xi^{1/2}\\
1 & 1
\end{pmatrix}
\begin{pmatrix}
E_{+,2n+1}(x,\xi)\\
E_{-,2n+1}(x,\xi)
\end{pmatrix}.
\]
Since the determinant of the $2\times 2$ matrix is $2\xi_*/\xi^{1/2} \neq 0$, solving for $E_{\pm,2n+1}(x,\xi)$, taking the limit \eqref{eq:SESL}, and then applying \eqref{eq:F1Fpm} and \eqref{eq:F4Fpm} proves \eqref{eq:EpmFpm} for {\em odd} dimensions.

Now, let $n$ be even. The equations \eqref{eq:E1Epm} and \eqref{eq:E2Epm} give the linear system
\[
\begin{pmatrix}
E_{1,n}(x,\xi_*)\\
E_{2,n}(x,\xi)
\end{pmatrix}
=
\frac{1}{2}\begin{pmatrix}
1+\xi_*/\xi^{1/2} & 1-\xi_*/\xi^{1/2}\\
E_{-,n+1}(x;\xi) & E_{+,n+1}(x;\xi)
\end{pmatrix}
\begin{pmatrix}
E_{+,n}(x,\xi)\\
E_{-,n}(x,\xi)
\end{pmatrix}.
\]
Since \eqref{eq:EpmFpm} is true for the odd dimension $n+1$, the determinant of the $2\times 2$ matrix becomes, in the soft-edge scaling limit,
\[
\Delta(s;\xi):=F_+(s;\xi) \bigl(1+\xi_*/\xi^{1/2}\bigr) - F_-(s;\xi) \bigl(1-\xi_*/\xi^{1/2}\bigr).
\]
By \eqref{eq:Fpm}, $\Delta(s;\xi)$ is zero if and only if
\[
\exp\biggl(-\int_{s}^\infty q(t;\xi)\,{\rm d}t\biggr) = \frac{\det\bigl(I-\xi^{1/2} V_{\Ai}\bigr)\big|_{L^2(s,\infty)}}{\det\bigl(I+\xi^{1/2} V_{\Ai}\bigr)\big|_{L^2(s,\infty)}} = \frac{1-\xi_*/\xi^{1/2}}{1+\xi_*/\xi^{1/2}} = \exp\bigl(-\artanh \xi^{1/2}\bigr).
\]
This is only possible for $0<\xi<1$ and then, by analyticity, at most in a discrete set of critical values of $s$. However, Figure~\ref{fig:2} indicates that such sets of $s$ exist for all $0<\xi<1$. Excluding those values of $s$, solving for $E_{\pm,n}(x,\xi)$ -- the linear systems is, for large $n$ in the soft-edge scaling limit, eventually solvable since the limit determinant is nonzero --, taking the limit~\eqref{eq:SESL}, and then applying \eqref{eq:F1Fpm} and \eqref{eq:F2Fpm} proves \eqref{eq:EpmFpm} also for {\em even} dimensions.
\end{proof}

\begin{Remark} Figure~\ref{fig:2} strongly suggests that the Ablowitz--Segur solution $q(s;\xi)$ of Painlev\'e~II has the total (improper) integral
\begin{equation}\label{eq:total}
\int_{-\infty}^\infty q(s;\xi) \,{\rm d}s = \artanh \xi^{1/2}, \qquad 0\leq \xi <1 .
\end{equation}
This is indeed the case, as demonstrated in \cite[equations~(27/28)]{MR2501035} by a Riemann--Hilbert analysis.
\end{Remark}

\subsection{Self-consistent expansion hypothesis} We generalize the scaling limits in Lemma~\ref{lem:EpmFpm} by embedding them into asymptotic expansions. In analogy to \cite[Section~4.2]{B25-1}, we hypothesize as follows.
\begin{HypothesisI} With $n_-:=n-1/2$, as $n\to \infty$, there holds the expansion\footnote{We are pretty confident that, with persistence and proper bookkeeping, based on the expansions of the Hermite and Laguerre polynomials of \cite[Part III]{B25-1}, one can prove Hypothesis I by extending the calculations in the work of Johnstone and Ma \cite{MR3025686,MR2888709} -- who established the case $\xi=1$ and $m=0$ for $\GOE$ and $\LOE$ using Pfaffian representations. Note that, by \eqref{eq:EgenRel1} and \eqref{eq:Epm}, expanding $E_{1,n}$ yields an expansion of $E_{\pm,n}$.}
\begin{gather}
E_{\pm,n}(x;\xi)|_{x=\mu_{n_-} + \sigma_{n_-} s}\nonumber\\
\qquad{} = F_\pm(s;\xi) + \sum_{j=1}^m G_{\pm,j}(s,\tau_{n_-};\xi)h_{n_-}^j + h_{n_-}^{m+1}O({\rm e}^{-s}),\qquad n\to\infty,\label{eq:Epmexpand}
\end{gather}
uniformly for $0\leq \xi \leq 1$ while $s$ stays bounded from below and, in the Laguerre case,~$\tau_{n_-}$~stays bounded away from zero -- the Gaussian case being encoded with $\tau_{n_-}=0$.

Here, the $G_{\pm,j}(s,\tau;\xi)$ are certain smooth functions and, preserving uniformity, the expansion can be repeatedly differentiated with respect to the variable $\xi$ and $s$.
\end{HypothesisI}

The shift in $n$ to $n_-$ in the expansion parameters is explained by $E_{2,n} (x;\xi)$ being the symmetrized product \eqref{eq:E2Epm}, mixing the dimensions $n$ and $n+1$ on the right-hand side. By observing
\[
n = \frac{n_- + (n+1)_-}{2},
\]
we see as in \cite[Section~4.2]{B25-1} that the hypothesized expansion \eqref{eq:Epmexpand}, if inserted to \eqref{eq:E2Epm}, consistently reclaims the form of the already {\em proven} expansion \eqref{eq:E2expand}. This is, if we re-expand quantities expressed in $(n+1)_-$ and $n_-$ to those expressed in $n$, we introduce terms in odd~powers of~\smash{$h_n^{1/2}$} which cancel by symmetry.

By spelling out the details of computing the symmetrized product \eqref{eq:E2Epm} of the self-consistent asymptotic expansions, we get equations relating $G_{2,j}$ with $F_\pm$, $G_{\pm,k}$, $k=1,2\dots,j$, -- such as for instance \cite[equation~(4.12)]{B25-1}:
\begin{equation}\label{eq:G21algebra}
G_{2,1} =F_+ G_{-,1} +G_{+,1}F_- + \frac{1}{2}\bigl(F_+ F_-'' - 2F_+'F_-' + F_+'' F_-\bigr).
\end{equation}
However, those equations are, in themselves, insufficient to determine the functions $G_{\pm,j}$.

\subsection{Linear form hypothesis} The determination of the functions $G_{\pm,j}$ becomes possible from a second structural
hypothesis. It~is important to note that, although this hypothesis allows us to compute the {\em assumed} polynomial coefficients and establish their independence of the dummy variable $\xi$ (at least for $j \leq m_*$), the fact that these computations can be carried out does not by itself imply that the hypothesis is true (for this range of $j$).

\begin{HypothesisII} The functions $G_{\pm,j}$ can be represented in the form
\begin{equation}\label{eq:Gpm}
G_{\pm,j}(s,\tau;\xi) = \sum_{k=1}^{2j} P_{\pm,j,k;\xi}(s,\tau) \cdot \partial_s^k F_{\pm}(s;\xi), \qquad P_{\pm,j,k;\xi} \in \Q[s,\tau].
\end{equation}
Note that we admit a possible dependency of the polynomial coefficients $P_{\pm,j,k;\xi}$ on $\xi$ here.
\end{HypothesisII}

The background of this hypothesis is that \eqref{eq:Fpm} implies, as in the proof of Theorem~\ref{thm:unitary},\footnote{\label{fn:algStruc}Note that this algebraic structure is lost for $F_1$, $F_4$, or the scaling limits of $E_{\even,n}$, $E_{\odd,n}$.}
\begin{equation*}
\frac{F_2^{(k)}}{F_2}, \frac{F_\pm^{(k)}}{F_\pm} \in \Q[s,q,q'].
\end{equation*}
Here, we use the notation of that proof, writing also $F_\pm^{(k)}=\partial_s^k F_\pm(s;\xi)$ for brevity.
We thus get, for instance,
\begin{gather}
\frac{F_\pm'}{F_\pm} = \frac{1}{2} q'^2-\frac{1}{2} q^4-\frac{s}{2} q^2\pm\frac{1}{2}q,\label{eq:Fpm2nd}\\
\frac{F_\pm''}{F_\pm} = -\frac{1}{2} q^4 q'^2-\frac{s}{2} q^2 q'^2\pm\frac{1}{2} q q'^2+\frac{1}{4} q'^4\pm\frac{1}{2}q'+\frac{s^2}{4} q^4 +\frac{1}{4}q^8 +\frac{s}{2} q^6\mp\frac{1}{2}q^5\mp\frac{s}{2} q^3-\frac{1}{4}q^2.\nonumber
\end{gather}
All those polynomial representations in $\Q[s,q,q']$ are {\em independent} of~$\xi$: the functional dependence on $\xi$ enters only indirectly via the boundary condition of the function~$q$ as $s\to\infty$. By the same arguments as in the proof of Theorem~\ref{thm:unitary}, we can apply literally the same calculations as in \cite[Section~4.2]{B25-1}: solving overdetermined linear systems over $\Q[s,\tau]$ that are obtained from comparing coefficients in $\Q[s,\tau,q,q'] =\Q[s,\tau][q,q']$ we get expressions for
the polynomials $P_{\pm,j,k;\xi}$ in terms of the coefficient polynomials $P_{2,l,m} \in \Q[s,\tau]$ of the unitary case -- expressions, which are by construction {\em independent} of $\xi$ and which happen to be, quite remarkably, also {\em independent} of the sign-index $\pm$.

For instance, combining \eqref{eq:G21algebra} with \eqref{eq:G2struct} and \eqref{eq:Gpm}
yields
\[
P_{2,11} \frac{F_2'}{F_2} + P_{2,12} \frac{F_2''}{F_2} = \sum_{\sigma = \pm} \biggl(P_{\sigma,11;\xi} \frac{F_\sigma'}{F_\sigma} + P_{\sigma,12;\xi} \frac{F_\sigma''}{F_\sigma}\biggr).
\]
Now, inserting the polynomial representations \eqref{eq:PII}, \eqref{eq:F22nd}, and \eqref{eq:Fpm2nd} gives, after comparing the coefficients of the 13 different $(q,q')$-monomials that enter the expressions,
\[
q, \, q^2,\, q^3,\, q^4,\, q^5,\, q^6,\, q^8,\, q',\, q'^2,\, q'^4,\, q q'^2,\, q^2q'^2,\, q^4q'^2,
\]
the $13\times 4$ linear system over $\Q[s,\tau]$ displayed in \cite[equation~(4.19)]{B25-1}. The unique solution is
\[
P_{\pm,11;\xi} = P_{2,11},\qquad P_{\pm,12;\xi} = 2P_{2,12}.
\]
Continuing in the same fashion, we get
\begin{alignat*}{3}
&P_{\pm,21;\xi} = P_{2,21}+2 P_{2,24} -\tfrac{1}{2} \partial_s^2 P_{2,11}+P_{2,12} - \tfrac{1}{4},\qquad&&
P_{\pm,22;\xi} = 2 P_{2,22} -\tfrac{1}{2} P_{2,11}^2-\partial_s^2 P_{2,12},&\\
&P_{\pm,23;\xi} = 4 P_{2,23}-2 P_{2,11} P_{2,12},\qquad&&
P_{\pm,24;\xi} = 8 P_{2,24} -2 P_{2,12}^2.&
\end{alignat*}
So, at least for as long as we keep continuing further on, that is up to at least $j\leq m_*$, we get
\begin{equation}\label{eq:P+}
P_{+,j,k;\xi}(s,\tau) = P_{-,j,k;\xi}(s,\tau) =: P_{1,j,k}(s,\tau) \qquad \text{(independent of $\xi$)}.
\end{equation}
The concrete instances up to $j=4$ are displayed in Table~\ref{tab:2}.

Finally, if we insert the asymptotic expansion \eqref{eq:Epmexpand} -- with the expansion terms $G_{\pm,j}$ given as the multinear form \eqref{eq:Gpm} with certain $\xi$-independent coefficient polynomials \eqref{eq:P+} -- into the relations \eqref{eq:E1Epm} and \eqref{eq:E4Epm}, we get by linearity and the limit relations \eqref{eq:F1Fpm} and \eqref{eq:F4Fpm}
the following (conditional) theorem. Regarding the shifts in $n$ we combine the cases $\beta=1,4$ into a~common notation by defining $n'$ as in~\eqref{eq:nprimeDef}~-- namely,
\[
n' =
\begin{cases}
n_- = n-\frac12, & \beta=1,\\
(2n+1)_- = 2n+\frac12 , & \beta=4.
\end{cases}
\]

\begin{Theorem}[conditional on Hypotheses I and II]\label{thm:orthogonal} There holds, for $\beta=1,4$ and any fixed non-negative integer $m \leq m_*$,\footnote{See footnote~\ref{fn:mstar}.}
\begin{gather*}
E_{\beta,n}(x;\xi) |_{x=\mu_{n'} + \sigma_{n'} s} = F_\beta(s;\xi) + \sum_{j=1}^m G_{\beta,j}(s,\tau_{n'};\xi)h_{n'}^j + h_{n'}^{m+1}O({\rm e}^{-s}),\qquad n\to\infty,
\end{gather*}
uniformly for $0\leq \xi \leq 1$ while $s$ stays bounded from below and, in the Laguerre case,~$\tau_{n'}$~stays bounded away from zero -- the Gaussian case being encoded with $\tau_{n'}=0$. Preserving uniformity, the expansion can be repeatedly differentiated with respect to the variables $s$ and $\xi$. The leading term $F_\beta(s,\xi)$ is given by \eqref{eq:Fbeta_Fpm}
and the correction terms are the multilinear forms
\begin{equation*}
G_{\beta,j}(s,\tau;\xi) = \sum_{k=1}^{2j} P_{\beta,j,k}(s,\tau)\cdot \partial_s^k F_\beta(s;\xi)
\end{equation*}
with certain polynomials $P_{\beta,j,k} \in \Q[s,\tau]$ that are independent of the dummy variable~$\xi$ and satisfy the duality relation $P_{1,j,k} = P_{4,j,k}$ -- the instances up to $j=4$ are displayed in Table~{\rm \ref{tab:2}}.
\end{Theorem}

Our substantiating evidence for that conditional theorem, and thus for the underlying structural Hypotheses I and II (which are, turning backwards, themselves implied by the statement of Theorem~\ref{thm:orthogonal} as the proof of Lemma~\ref{lem:EpmFpm} shows), was already collected in Section~\ref{subsect:evidence} -- see also Figure~\ref{fig:evidence}.

\appendix

\section{Parametrization of the matrix ensembles}\label{app:ensembles}
For convenience, we recall from \cite{B25-1,B25-3} our choice of parametrization. The joint probability density of the (unordered) levels $x_1,\dots,x_n$ of an $n$-dimensional $\beta$-ensemble (with $\beta=1,2,4$ encoding
 the orthogonal, unitary, or symplectic cases) takes the form
\begin{equation*}
P_\beta(x_1,\dots,x_n) \propto |\Delta(x_1,\dots,x_n)|^\beta \prod_{j=1}^n w_\beta(x_j),
\end{equation*}
where $\Delta(x_1,\dots,x_n) = \prod_{j>k} (x_j - x_k)$ denotes the Vandermonde determinant and the weight functions $w_\beta(\xi)$ are given, for the Gaussian ensembles, by the Hermite weights
\begin{subequations}\label{eq:alphabeta}
\begin{equation}
w_\beta(x) = {\rm e}^{-c_\beta x^2}, \qquad x\in\R
\end{equation}
and, for the Laguerre ensembles with $\alpha>-1$, by the Laguerre weights
\begin{equation}
w_\beta(x) = x^\alpha {\rm e}^{-c_\beta x}, \qquad x>0.
\end{equation}
However, we replace the {\em Laguerre parameter} $\alpha$ with the more convenient {\em Wishart parameter}
\begin{equation}\label{eq:palpha}
p= n - 1 + \frac{2 (\alpha + 1)}{\beta} > n-1.
\end{equation}
If $p$ is an integer, the levels of the Laguerre ensemble are distributed as the eigenvalues of the $n$-variate Wishart distribution with $p$ degrees of freedom. Since the Wishart distribution exhibits a~$p\leftrightarrow n$ symmetry \cite[Section~1]{B25-1}, the formulas for the Laguerre ensemble can be put to a form reflecting that symmetry, which results in considerable simplifications. Consequently, we write $\GbE_n$ for the Gaussian ensembles and $\LbE_{n,p}$ for the Laguerre ensembles. Fixing the scales of the levels by
\begin{equation}
c_\beta = \begin{cases}
1/2, & \beta = 1,\\
1, & \beta = 2,4,
\end{cases}
\end{equation}
\end{subequations}
we facilitate the Forrester--Rains interrelations \eqref{eq:ForresterRains}.

\section{Gap-probability generating functions of point processes}\label{app:PGF}

We follow \cite{MR854102} and define a point process as an integer-valued random measure $N$ (the counting measure) on a Polish space $\Omega$. For any Borel measurable $A\subset \Omega$, $N(A)$ is thus a random variable taking values in the nonnegative integers. The associated probability-generating function (PGF) can be written in the form\footnote{We take the dummy variable $\xi$ of the PGF in this way instead of the more obvious choice $\E\bigl(x^{N(A)}\bigr)$ to match with the common definition in the random matrix literature (see, e.g., \cite[Section~8.1]{MR2641363}).}
\begin{equation*}
E(A;\xi):= \E\bigl((1-\xi)^{N(A)}\bigr) = \sum_{k=0}^\infty \prob (N(A)=k ) \cdot (1-\xi)^k,
\end{equation*}
where $\xi$ is often addressed as the formal, indeterminate, or dummy variable. In random matrix theory, the probability $\prob(N(A)=0)$ -- that there is no point in $A$ -- is called the gap probability. Accordingly, we refer to $\prob(N(A)=k)$ as the $k$-level gap probability and call~$E(A;\xi)$ the gap-probability generating function associated with the set $A$.

For $0\leq\xi \leq 1$, the value of the generating function $E(A;\xi)$ has a simple probabilistic interpretation itself (as is well known in the theory of renewal and branching processes): if we keep each point of the point process independently with probability $\xi$ and delete it otherwise, such a $\xi$-(Bernoulli-)thinning yields (cf., e.g., \cite[Table 1, p.~16]{MR854102})
\[
\prob(\text{no point of the $\xi$-thinning belongs to $A$}) = E(A;\xi).
\]
By Taylor's formula, the $k$-level gap probability can be reconstructed by evaluating derivatives at $\xi=1$,
\begin{equation*}
\prob(N(A)=k) = \frac{(-1)^k}{k!} \frac{{\rm d}^k}{{\rm d}\xi^k} E(A;\xi)\biggr|_{\xi=1}.
\end{equation*}
On the other end, at $\xi=0$, we get
\begin{equation*}
\E(N(A)) = \sum_{k=0}^\infty k \cdot \prob(N(A)=k) = - \frac{{\rm d}}{{\rm d}\xi} E(A;\xi)\biggr|_{\xi=0}.
\end{equation*}
We now specialize to continuous finite simple point processes on a closed subset $\Omega \subset \R^d$, realized by $n$ a.s. distinct points $x_1,\dots,x_n \in \Omega$ that have a continuous and symmetric joint probability density function $P(x_1,\dots,x_n)$. Upon writing\footnote{Here, $\llbracket\cdot \rrbracket$ denotes the Iverson bracket: $\llbracket \mathcal P \rrbracket=1$ if the proposition $\mathcal P$ is true, and $\llbracket \mathcal P \rrbracket=0$ otherwise.}
\[
N(A) = \sum_{j=1}^n  \llbracket x_j \in A \rrbracket, \qquad (1-\xi)^{N(A)} = \prod_{j=1}^n (1-\xi)^{\llbracket x_j \in A \rrbracket},
\]
and noting that
\[
(1-\xi)^{\llbracket x_j \in A \rrbracket} = 1-\xi \llbracket x_j \in A \rrbracket,
\]
we immediately get an alternative form of the generating function (cf., e.g., \cite[Proposition~8.1.2]{MR2641363})
\begin{subequations}\label{eq:EAint}
\begin{align}
E(A;\xi) &= \int_\Omega\cdots \int_\Omega P(x_1,\dots,x_n) \prod_{j=1}^n \bigl(1-\xi \llbracket x_j \in A \rrbracket\bigr) \,{\rm d}x_1\cdots {\rm d}x_n\nonumber\\
&=1 + \sum_{k=1}^n \frac{(-\xi)^k}{k!} \int_A\cdots \int_A R_k (x_1,\dots,x_k)\,{\rm d}x_1\cdots {\rm d}x_k,
\end{align}
where the $k$-point functions $R_k$ are defined by\footnote{Here, we use Knuth's notation $n^{\underline{k}} = n(n-1)\cdots(n-k+1)$ for the falling factorial power.}
\begin{equation}
R_k (x_1,\dots,x_k) := n^{\underline{k}} \int_\Omega\cdots \int_\Omega P(x_1,\dots,x_k,x_{k+1},\dots,x_n)\,{\rm d}x_{k+1}\cdots {\rm d}x_n.
\end{equation}
\end{subequations}
If the point process is also determinantal with correlation kernel $K(x,y)$, then the $k$-point functions are given by
\[
R_k(x_1,\dots,x_k) = \det_{\mu,\nu=1}^k K(x_\mu,x_\nu),\qquad k=1,2,\dots,n,
\]
and the polynomial expansion \eqref{eq:EAint} can briefly be recast as a Fredholm determinant (cf., e.g., \cite[Section~4.2]{MR2760897}, \cite[Section~5.4]{MR1677884} and \cite[Section~9.1]{MR2641363}) of the induced integral operator $K$ acting on~$L^2(A)$,
\begin{equation}\label{eq:EAdet}
E(A;\xi) = \det(I-\xi K)|_{L^2(A)}.
\end{equation}

\subsection*{Acknowledgements} Special thanks to Peter Forrester for suggesting that I generalize my work \cite{B25-1} to generating functions and, more generally, for conversations during a visit to Melbourne he hosted in March 2024 -- a period in which the results of this paper took shape.

\pdfbookmark[1]{References}{ref}
\LastPageEnding

\end{document}